\documentclass{article}
\usepackage{amsmath,amsfonts,amssymb,theorem,verbatim}
\newtheorem{theorem}{Theorem}%[section]
\newtheorem{lemma}[theorem]{Lemma}%[section]

\newtheorem{prop}[theorem]{Proposition}
\def\BEN{\begin{enumerate}}  \def\BI{\begin{itemize}}
\def\EEN{\end{enumerate}}   \def\EI{\end{itemize}}
    \def\nn{\nonumber}

\def\beq{\begin{eqnarray}} \def\eeq{\end{eqnarray}}

\def\eqn#1{\begin{equation}#1\end{equation}}

\def\al*#1{\begin{align*}#1\end{align*}}

\def\ga*#1{\begin{gather*}#1\end{gather*}}

\def\alat*#1#2{\begin{alignat*}{#1}#2\end{alignat*}}
\def\bea{\begin{eqnarray*}}
\def\eea{\end{eqnarray*}}
\def\ml*#1{\begin{multline*}#1\end{multline*}}

 \def\mbf{\mathbf} 
\def\mc{\mathcal} \def\unl{\underline} \def\ovl{\overline}
\def\tld{\tilde}

\def\P{{\mathbb P}} \def\le{\left} \def\ri{\right} \def\i{\infty}
   \def\R{{\mathbb R}}

\def\FF{\mathcal{F}}  
   
\def\maxe{\overline{\epsilon}} \def\mine{\underline{\epsilon}}

\def\te#1{\mathrm{e}^{#1}}   
\def\WT{\widetilde}

  \def\im{\item} 
\def\I{\int}     \def\a{\alpha} 
  \def\d{\delta} \def\z{\zeta}  \def\th{\theta}

\def\e{\epsilon} \def\k{\kappa}  \def\m{\mu} 
\def\x{\xi}  \def\nn{\nonumber}   \def\s{\sigma}
\def\t{\tau}   \def\f{\varphi} \def\x{\xi} \def\ps{\psi}
 \def\w{\omega} \def\q{\qquad} \def\D{\Delta}
\def\F{\Phi}  \def\L{\Lambda} \def\O{\Omega} 
  \def\td{\text{\rm d}}
\numberwithin{equation}{section}

 \newtheorem{As}{Assumption}
{\theorembodyfont{\normalfont}

}

\newcommand{\proof}{{\it Proof\ }}

\newcommand{\exit}{{\mbox{\, \vspace{3mm}}} \hfill\mbox{$\square$}}

\begin{document}
\date{February 2008}
\title{On an explicit Skorokhod embedding for spectrally
negative L\'{e}vy processes}
\author{Jan Ob{\l}\'{o}j\thanks{e-mail:
        \texttt{jobloj@imperial.ac.uk}; web:
        \texttt{www.imperial.ac.uk/people/j.obloj/}
        \newline Research supported by a Marie Curie Intra-European
Fellowship within the $6^{th}$ European Community Framework Programme}\\
        Dept.\ of Mathematics\\
        Imperial College London\\
        London SW7 2AZ, UK
\and Martijn Pistorius\thanks{e-mail: \texttt{Martijn.Pistorius@kcl.ac.uk};
web: \texttt{www.mth.kcl.ac.uk/staff/m\_pistorius.html}
\newline Research supported by EPSRC grant EP/D039053/1}\\
Dept.\ of Mathematics\\
King's College London\\
London WC2R 2LS, UK
}
\maketitle
\begin{abstract}
We present an explicit solution to the Skorokhod embedding problem
for spectrally negative L\'evy processes. Given a process $X$ and
a target measure $\mu$ satisfying an explicit admissibility
condition we define functions $\f_\pm$ such that the stopping time
$T = \inf\{t>0: X_t \in \{-\f_-(L_t), \f_+(L_t)\}\}$ induces
$X_T\sim \mu$, where $(L_t)$ is the local time in zero of $X$. 
We also treat versions of $T$ which take into
account the sign of the excursion straddling time $t$. We prove
that our stopping times are minimal and we describe criteria under
which they are integrable. We compare our solution with the one
proposed by Bertoin and Le Jan \cite{blj}. In particular, we
compute explicitly the quantities introduced in \cite{blj} in our
setup.

Our method relies on some new explicit calculations relating scale
functions and the It\^o excursion measure of $X$. More precisely,
we compute the joint law of the maximum and minimum of an
excursion away from 0 in terms of the scale function.

\end{abstract}

\section{Introduction}
The \emph{Skorokhod embedding problem} was first introduced and solved by Skorokhod \cite{skoro}, where it
served to realize a random walk as a Brownian motion stopped at a sequence of stopping times. Since
then, it remains an active field of study and the original problem has been generalized in a number
of ways and has known many different solutions.
We refer to Ob\l\'oj \cite{genealogia} for a comprehensive survey paper.

The embedding problem can be phrased in the following general manner: given a stochastic process
$(X_t)$ and a measure $\mu$ on its state space, find a stopping time $T$ which embeds the measure: $X_T\sim
\mu$. To make the problem interesting one requires $T$ to be \emph{small} in some sense.

When $(X_t)$ is a continuous martingale and $\mu$ is centered, one typically asks that $(X_{T\land t})$ is a uniformly integrable martingale.
When $\mu$ has finite second moment this is equivalent to the expectation of the quadratic variation
being finite, $E[\langle X\rangle_T]<\i$, which for Brownian motion reads simply $E[T]<\i$. A more general condition requires $T$
to be minimal, that is if $S$ is a stopping time with $S\leq T$ and $X_S\sim X_T$ then $S=T$. Monroe \cite{MR49:8096} showed that the two conditions
are equivalent when $X$ is a continuous local martingale and $\mu$ is a centered probability measure.
In contrast, Cox and Ob\l\'oj \cite{cox_obloj_fractal} showed that for discontinuous processes, even in the simplest
case of a symmetric random walk, this no longer holds true. In fact, these authors showed that the set of measures which can be embedded in a uniformly integrable way may be a complex fractal subset of the
set of measures which can be embedded using minimal stopping times.

In this paper we solve the Skorokhod embedding problem for spectrally negative L\'evy processes.
Our solution is based on the general framework developed by Ob\l\'oj \cite{ob_gen}, and recently
used in Pistorius \cite{mp-sem07}. Here, we have to extend the setup
to account for the presence of jumps.
In order to do so we carry out excursion theoretical computations which have an interest in their own. More precisely, given a target measure $\mu$ from a certain class, we find
functions $\f_\pm$ such that the stopping time
\begin{eqnarray}\label{eq:T}
T_{\f_\pm} &=& \inf\{t>0: X_t \in \{-\f_-(L_t), \f_+(L_t)\}\},
%\\ \tld T_{\f_\pm}&=&\inf\{t>0: X_t=\varphi_+(L_t)\textrm{ or }X_t=-\varphi_-(L_t)\textrm{ and }sgn(e_{L_t})=-1\},\label{eq:Tbis}
\end{eqnarray}
where $L_t$ is the local time at zero, satisfies $X_{T_{\f_\pm}}\sim\mu$. We also look at versions
of \eqref{eq:T} which take into account the sign of the excursion straddling
time $t$.
We compute the joint law of the maximum and minimum of an excursion
away from zero in terms of the scale function, and this yields the functions $\f_{\pm}$ explicitly in terms of $\mu$ and the characteristics of $X$.
Finally we prove that our stopping times are minimal
-- and this in spite of possible waiting for a 'comeback' after an undershoot.
To the best of our knowledge, so far only one solution to our embedding problem exists,
proposed by Bertoin and Le Jan \cite{blj}. In \cite{blj}
a general solution is presented and then the case of symmetric L\'evy processes
is treated in detail. We complement this by providing explicit computations
for the case of spectrally negative L\'evy processes.

In \cite{blj} stopping times are constructed
using local times at all levels simultaneously, while our construction uses only the local
time in zero. The difference in complexity is best seen when $X_t=B_t$ is a Brownian motion,
comparing the formulae in \cite{blj} with the solution of Vallois \cite{val}, to which our
solution reduces in this simplest setup.
Naturally, there is a price we have to pay for having a simple explicit construction, namely
we can only embed measures satisfying a certain admissibility criterion. We discuss
this in detail and state the criterion in terms of the scale function.

The paper is organised in the following manner. In Section \ref{sec:levy} we describe the class of L\'evy processes we will be working with and the associated objects such as the scale functions and local times. In Section \ref{sec:skoro} we state our embedding theorems
and discuss the minimality of stopping times and the class of admissible measures. The subsequent section contains excursion theoretical results, namely the computations which are the main ingredient for our formulae. Finally in Section \ref{sec:proofs} we prove the embedding results.

\section{Preliminaries}\label{sec:levy}
Let $X=(X_t, t\ge0)$ be a spectrally negative L\'{e}vy process defined on
a filtered probability space $(\O,\FF,{\mbf F}=\{\FF_t\}_{t\ge0},P)$.
Here the filtration $\mbf F$ is the completion of the standard filtration
generated by $X$. To avoid trivialities,
we exclude the case that $X$ has monotone paths. Since the  jumps
of $X$ are all non-positive, the moment generating function
$E[\te{\th X_t}]$ exists for all $\th\ge0$ and is given by
$\ps(\th) = t^{-1} \log E[\te{\th X_t}]$ for some function
$\ps(\th)$. The function $\psi$ is well defined at least on
the positive half-axis where it is strictly convex with the
property that $\lim_{\th\to\i}\ps(\th)=+\i$. Moreover, $\psi$
is strictly increasing on $[\F(0),\i)$,
where $\F(0)$ is the largest root of $\ps(\th)=0$. We shall denote
the right-inverse function of $\ps$ by $\F:[0,\i)\to[\F(0),\i)$.
Note that $\F(0)>0$ if and only if $X$ drifts to $-\i$ (see Bertoin \cite[Cor.~VII.2]{bert}).

The continuous martingale component of $X$ is a Brownian
motion with variance $\s^2$, called the Gaussian coefficient
 of $X$, which can be recovered from $\psi$ by $\s^2 =
\lim_{\th\to\i} 2\psi(\th)/\th^2$. In what follows we make the
following assumption on the process $X$:
\begin{As}
\label{as:levy}
$X$ is a spectrally negative L\'{e}vy process that has unbounded
variation and does not drift to $-\i$.
\end{As}

The hitting time of a set $\Gamma\subset \R$ is denoted $H_\Gamma=\inf\{t: X_t\in \Gamma\}$.
We write $H_\eta$ and $H_{\eta,\delta}$ respectively for $H_{\{\eta\}}$ and $H_{\{\eta,\delta\}}$.
An important role in the fluctuation theory of spectrally
negative L\'{e}vy processes is played by the so-called {\it $q$-scale functions}
$W^{(q)}: \R \to [0,\i)$, $q\ge0$, that are zero on the negative half-axis and
continuous and increasing on $[0,\i)$ with Laplace transforms
\eqn{\label{eq:defW}
\I_0^\i \te{-\th x} W^{(q)} (x)  \td x = (\ps(\th) - q)^{-1},\q\th > \F(q).
}
See e.g. Bingham \cite{bingham}, Bertoin \cite[Thm.~VII.8]{bert}
or Kyprianou \cite[Thm 8.1]{Kbook}
for proofs of the existence of $W^{(q)}$.
Bertoin \cite{bertoin:exitergo}
has shown that for every $x\ge0$,
the mapping $q\mapsto W^{(q)} (x)$ can be analytically
extended to the complex plane by the identity
\eqn{\label{eq:expW}
W^{(q)} (x) = \sum_{k\ge0}q^kW^{\star (k+1)} (x),
}
where $W^{\star k}$ denotes the $k$-th convolution power of
$W=W^{(0)}$.

The scale function $W$ is closely linked to the law of $X_{H_{\R \setminus [a,b]}}$:
\begin{equation}\label{twosided}
P_x[X_{H_{\R \setminus [a,b]}} = b] = \frac{W(x-a)}{W(b-a)},\qquad a<x<b,
\end{equation}
where $P_x[\cdot] = P[\cdot|X_0=x]$. Moreover, from Bertoin
\cite[Cor. 1]{bertoin:exitergo} it follows that
\begin{equation}\label{twosidedexp}
E_x[H_{\R \setminus [a,b]}] = \frac{W(x-a)}{W(b-a)}\ovl W(b-a) - \ovl W(x-a),
\qquad a<x<b,
\end{equation}
where $\ovl W(x) = \int_0^x W(y)\td y$. 
We note that the Gaussian
coefficient can be recovered from $W$ as $W'(0):=W'(0+)=2/\s^2$
(see \cite[Lemma 1]{mp-sem05}).

The $q-$potential measure $U^q(\td x) =
\int_0^\i\te{-qt}\P(X_t\in\td x)\td t$ of $X$ is absolutely continuous
with density $u^q$ related to the $q-$scale function $W^{(q)}$ by
\begin{equation}\label{eq:uq}
u^q(x) = \F'(q)\te{-\F(q)x} - W^{(q)}(-x), \quad q>0
\end{equation}
(see Bingham \cite{bingham} or Pistorius \cite{mp-sem05}). 
Since the measure $U^q$ is absolutely continuous with bounded density $u^q$
and $0$ is regular for $\{0\}$ if $X$ has unbounded variation
(\cite[Cor VII.5]{bert}), the limit, a.s.\ and in $L^2(P)$,
$$
L^y_t = \lim_{\e\downarrow 0} \frac{1}{2\e}\I_0^t \mathbf 1_{\{|X_s-y| < \e\}}\td s
$$
exists for every $y\in\R$ and $t\ge0$ 
and $t\mapsto L^y_t$ is continuous a.s. 
(cf.\ Bertoin \cite[Thm. II.16, Thm V.1, Prop V.2]{bert}), so that, in particular, $E_x [L^y_t]<\i$.
The process $(L^y_t, t\ge 0)$ is called the local time of $X$ at $y$. 
The expectation
of the Laplace-Stieltjes transform in time of the local time
$L^y_t$ is related to the potential density via
\begin{equation}\label{eq:uqy}
E\le[\I_0^\i\te{-qt}\td L^y_t\ri] = u^q(y).
\end{equation}
Denote by $\t = (\t_\ell, \ell \ge 0)$ the right-continuous 
inverse of $L = L^0$,
$$
\t_\ell = \inf\{t>0: L_t > \ell\}.
$$

\section{The Skorokhod embedding problem}\label{sec:skoro}

Henceforth we consider  the Skorokhod embedding problem for a measure
$\mu$ on $\R \setminus \{0\}$ under the condition of minimality on a solution $T$ (i.e.\
$R=T$ for any stopping time $R\leq T$ with $X_R\sim X_T$):
\begin{equation*}
(S) \quad
\begin{array}{l}
\text{Find a {minimal} almost surely finite $\mbf F$-stopping
time $T$} \\ \text{such that $X_T\sim\mu$}.
\end{array}
\end{equation*}
We will present two solutions to $(S)$:
one which works only in the presence of a positive
Gaussian component and a general one, which is however less explicit.
The two coincide and simplify for measures concentrated on $\R_+$. The explicit formulae in our solutions are a consequence of the excursion theoretical computations presented in Section
\ref{sec:exc_comp}.
\bigskip
\bigskip

\subsection{Solution in the presence of a Gaussian component}\label{sec:sol1}
Assume $\s^2>0$.  In this case our solution extends the ideas of Ob\l\'oj \cite{ob_gen}
to the discontinuous setup and in particular it simplifies to Vallois's solution \cite{val}
when $X$ is a Brownian motion.
We follow the approach of Cox, Hobson and Ob\l\'oj \cite{cho} to account for measures with atoms.
We note also that an atom in zero can easily be treated (see
\cite{val,cho}).

Let $a_\m<0<b_\m$ be the infimum and supremum of the support of $\m$
respectively and denote by $F_\mu(x)=\mu((-\i,x])$ the
cumulative distribution function of $\mu$, $F^{-1}_\mu$
its right-continuous inverse and let $a_*=F_\mu(0)$. We impose the
following \emph{admissibility} criterion on the measure $\mu$:
\begin{equation}\label{eq:assum}
\I_0^\i W(s)\m(\td s) =W'(0)\I_{-\i}^0\frac{W(-s)\m(\td s)}{W'(-s)}.
\end{equation}

Define $\alpha:[a_*,1]\to [0,a_*]$ via
\begin{equation}\label{eq:defalpha}
\I_{a_*}^a W(F^{-1}_\mu(s))\td s=W'(0)\I_{\alpha(a)}^{a_*}
\frac{W(-F^{-1}_\mu(s))}{W'(-F^{-1}_\mu(s))}\td s.
\end{equation}
Note that $\alpha$ is a strictly decreasing, absolutely continuous
function with
$\a(a_*)=a_*$ and $\a(1)=0$. Define $\x=\x_\mu$ via
\[ \x(a) = \int_{a_*}^{a} \frac{ W(F^{-1}_\mu(s)) }{\alpha(s) + (1-s)} \td s
\hspace{38mm} a_* \leq a \leq 1 \]
and
\[ \x(a) = \int_{a}^{a_*}
\frac{-W'(0)W(-F^{-1}_\mu(s))\td s}{W'(-F^{-1}_\mu(s))
(s + 1- \alpha^{-1}(s))}
\hspace{10mm} 0 \leq a \leq a_*  , \]
and put $\tld \psi_\mu(x)=\x(F_\mu(x))$.\footnote{A simplified expression for $\tld \psi$ in the case when $\mu$ has no atoms is given in Section \ref{sec:proofs}.}
Note that $\tld\psi$ is an increasing function on $\R$. We denote by $\tld\f$
the right-continuous inverse of $\tld\psi$ and write $\f_\pm(l)=\pm\tld\f(\pm l)$.

\bigskip

\begin{theorem}\label{thm:main}
Let $\s>0$ and suppose \eqref{eq:assum} holds. Then
\begin{equation}\label{eq:Tbis}
    \tld T_{\f_\pm}=\inf\{t>0: X_t=\varphi_+(L_t)\textrm{ or }X_t=-\varphi_-(L_t)
\textrm{ and } X_{|_{(\t_{L_t-},t)}} < 0\}
\end{equation}
solves $(S)$.
If $E[X_1^2]<\i$, $X$ drifts to $+\i$ and
\begin{equation}
\label{eq:cond_et}
\int_0^\i y W(y)\mu(\td y) - \int_{-\i}^0 y\frac{W(-y)}{W'(-y)}\m(\td y)<\i
\end{equation}
then $E[\tld T_{\f_\pm}]<\i$.
\end{theorem}

\subsection{Solution for $\sigma\geq 0$ and measures with a density}
Let $\s^2\geq 0$ be arbitrary.
Note that when there is no Gaussian component ($\s=0$) it is not possible to
separate the excursions into positive and
negative ones: every excursion starts positive and either
stays always positive or becomes negative and then ends.
We will restrain ourselves to probability measures $\mu$
with a positive density function $f_\mu$ on $(a_\mu,b_\mu)$, where $a_\mu$ and $b_\mu$ are respectively
the lower and the upper bound of the support of $\mu$.
Let $g:[0,b_\mu)\to(a_\mu,0]$ be a continuous decreasing function given
via
\begin{equation}
\label{eq:g_diff} 
\frac{\td g}{\td y}(y) = - \frac{W(y
- g(y)) - W(-g(y))}{W(-g(y))} \frac{f_\mu(y)}{f_\mu(g(y))}, \quad y\in(0,b_\mu),
\end{equation}
with $g(0) = 0$. We impose the following admissibility assumption on $\mu$
\begin{equation}\label{eq:assum_gen}
g\textrm{ is well defined and finite with } g(x)\xrightarrow[x\to b_\mu]{}a_\mu.
\end{equation}
Let $h$ denote the right-continuous inverse of $g$ and set for $y,z \ge 0$, $\ovl\mu(s) = \mu([s,\i))$,
\begin{equation}\label{eq:f0def}
\begin{split}
\psi_+(y) =& \I_0^y \frac{W(s)\mu(\td s)}{(1+\ovl\mu(s) - \ovl\mu(g(s)))}\\
\psi_-(z) =& \I_{-z}^0 \frac{W(h(s))W(-s)\mu(\td
s)}{\left(W(h(s)-s)-W(-s)\right) \left(1+\ovl\mu(h(s)) -
\ovl\mu(s)\right)}.
\end{split}
\end{equation}
\begin{theorem}\label{thm:main0}
Suppose $\mu$ is a probability measure on $\R$ with a positive
density $f_\mu$ on $(a_\mu,b_\mu)$ and which satisfies \eqref{eq:assum_gen}.
Then
\begin{equation}\label{eq:Tgen}
T_{\f_\pm} = \inf\{t>0: X_t \in \{-\f_-(L_t), \f_+(L_t)\}\},
\end{equation}
solves $(S)$, where $\f_+$, $\f_-$ are the right-continuous inverses of  respectively
$\psi_+$ and $\psi_-$ from \eqref{eq:f0def}.
\end{theorem}
Note that the functions $\f_\pm$ in \eqref{eq:Tbis} and \eqref{eq:Tgen} are different. It should
be clear which functions we mean depending on which stopping time is discussed.
\subsection{Solution for measures concentrated on $(0,\i)$}\label{sec:sol3}
We specialize now to the case where the target measure $\mu$ is a probability
measure on $(0,\i)$, in which case the solution further simplifies.
Set $T_\mu$ equal to the stopping time
\begin{equation}\label{eq:Tpos}
T_{\mu}=\inf\{t>0: X_t\geq \f_\mu(L_t)\}
\end{equation}
where $\f_\mu$ denotes the right-continuous inverse of
the map $\psi_\mu: [0,\i) \to [0,\i]$ that is defined by
\begin{equation}
\label{eq:fpos}
\psi_\mu(y)=-\int_{[0,y)}W(s)\td\big(\log\ovl\mu(s)\big).
\end{equation}

\begin{theorem}
\label{thm:pos}
Suppose $\m(0,\i)=1$. Then $T_\mu$ solves $(S)$ and
$\sup_{t\leq T_{\mu}} X_t=X_{T_\mu}$.
\end{theorem}

\bigskip

\subsection{On minimality of stopping times and admissibility of target measures}

We want to stress the minimality property of the stopping
times in Theorems \ref{thm:main} and \ref{thm:main0}.
It may seem surprising as the following example demonstrates.
Consider the first hitting time $H_\Gamma$ of a region
$\Gamma=(-\i,a]\cup [b,\i)$ for some $a<0<b$, which embeds a
distribution $\mu$ that has an atom in $b$ and the rest of the
mass in $(-\i,a]$. Now if we try to develop an embedding of $\mu$
with Theorems \ref{thm:main} or \ref{thm:main0} it seems that we
may very well stop later (also in the local time scale) than $H_\Gamma$.
The answer to this apparent paradox is that the measure $\mu$ can not
be treated within our framework.
The admissibility assumptions \eqref{eq:assum} and \eqref{eq:assum_gen} are crucial and they determine the set of measures that can be treated with our methodology. The restriction
of the set of admissible measures is a natural price to pay
for having a simple explicit form of
the stopping time that involves only the local time at zero.

Condition \eqref{eq:assum} requires that $\mu$ is \emph{centered} relative
to some density on $\R$ expressed in terms of the scale function. In particular when $X_t=B_t$ is a Brownian motion $W(s)=2s\mathbf{1}_{s\geq
0}$ and \eqref{eq:assum} simplifies to  $\int_{-\i}^0 |x|\mu(\td x)=\int_0^{\i} x\mu(\td x)$, which is a necessary and sufficient condition for $\f_\pm$ to be well defined and finite on $\R_+$. In general when $X$ is recurrent any measure on $\R$ can be embedded in $X$. An extension of our construction to
arbitrary measures would involve explosion of one of the functions $\f_\pm$ and the resulting
stopping times may not be minimal (as example described above illustrates).

When $X$ drifts to infinity it is not possible to embed all
measures on $\R$ in $X$. More precisely, Rost's balayage condition \cite{rost} states that there
exists an embedding of $\mu$ in $X$ if and only if
\begin{equation}\label{eq:rost}
E_0 \le[\int_0^\infty f(X_t)\td t\ri]\geq E_\mu\le[\int_0^\i f(X_t)\td t\ri],\quad \textrm{for all }f\geq 0.
\end{equation}
Naturally this is equivalent to a restriction on the potential
density of $X$, which we can rewrite using \eqref{eq:uq} as
$\int_\R (W(-y) - W(x-y))\mu(\td x) \leq 0$ a.e. Our embedding
works for measures which satisfy \eqref{eq:assum}, which is thus a
subclass of all measures which satisfy Rost's condition (note that
it may not be easy to show this directly).
\\
As an example, consider Brownian motion with drift $X_t=B_t+\delta t$, $\delta>0$. Then $W(s)=(1-\exp(-2\delta s))/\delta$, $s\geq 0$, and \eqref{eq:assum} simplifies to $\int_{\R}(1-\exp(-2\delta x))\mu(\td x)=:m=0$, while a necessary
and sufficient condition \eqref{eq:rost} on $\mu$ for existence of an embedding in $X$ is
$m\ge0$ (see Ob\l\'oj \cite[Sec.~9]{genealogia}). Note that if the last integral is equal
 to $m>0$ then we can still embed $\mu$ by using our
construction for the shifted process. More precisely, let $c=-\ln(1-m)/2\delta>0$ and $\tld \mu(\td
u)=\mu(\td(u+c))$. Then $\tld\mu$ satisfies \eqref{eq:assum} and embedding $\tld \mu$ in the
process $\tld X_t:= X_{t+H_c}-c$ we embed $\mu$ in $X$.\footnote{Naturally, similar arguments can be used to embed measures which violate \eqref{eq:assum} for a recurrent $X$, but again the resulting stopping times need not be minimal.}

For a continuous process
$X$ with $\sigma>0$ Theorems \ref{thm:main} and \ref{thm:main0} coincide.
However, for a discontinuous process $X$ with $\s>0$ the sets
of measures that can be embedded using Theorem \ref{thm:main}
and Theorem \ref{thm:main0} are different. We will come back to
this discussion when presenting examples in Section \ref{sec:exam}.

The key ingredient for the proof of minimality of our
stopping times is the observation that
they minimise the expectation of the local time at zero
among all stopping times which embed a
given law. This is closely related to the work of Bertoin
and Le Jan \cite{blj} and is discussed in the subsequent section.

\subsection{On the solution of Bertoin and Le Jan \cite{blj}}

Bertoin and Le Jan \cite{blj} presented a general solution to the
Skorokhod embedding problem. We explore now briefly their solution
in the context of spectrally negative L\'evy processes.

For a given probability measure $\mu$, Bertoin and Le Jan
\cite{blj} defined the function\footnote{Note that our $V_\mu$
corresponds to $\hat{V}_\mu$ in \cite{blj}.}
\begin{equation}\label{eq:VBLJ}
    V_\mu(y)=\int_{\R} v(x,y)\mu(\td x),\quad \textrm{and}
\quad \lambda_\mu=\sup_y V_\mu(y),
\end{equation}
where $v(x,y) = E_x[L^y_{H_0}]$. Bertoin and
Le Jan \cite{blj} proved that when $\lambda_\mu<\i$ the stopping time
\begin{equation}\label{eq:TBLJ}
    T_{BLJ}^\mu=\inf\left\{t:\int_{\R}\frac{\lambda_\mu L^x_t}{\lambda_\mu-V_\mu(x)}\mu(\td x)>L^0_t\right\}
\end{equation}
solves the Skorokhod embedding problem,
i.e.~$X_{T_{BLJ}^\mu}\sim \mu$, when $X$ is recurrent.
They also proved that among all solutions to the
Skorokhod embedding problem $T_{BLJ}^\mu$ minimizes the
expected value of additive functionals and that it holds that
$E[L^0_{T_{BLJ}^\mu}]=\lambda_\m$. We recover this bound in our setting:

\begin{prop}\label{prop:min}
(i) For any a.s. finite stopping time $S$ with $X_S\sim \mu$,
$$
E[ L_S]\geq \lambda_\mu\geq \int_0^\i W(x)\mu(\td x).
$$

(ii) If $X$ drifts to $+\i$, then $\lambda_\mu<\i$. 

(iii) Assume $X$ oscillates and that $\sigma^2>0$ and let $\mu$ be a probability measure
satisfying \eqref{eq:assum}. If $\int_0^\i x\mu(\td x)<\i$ then $\lambda_\mu<\i$. 
If, in addition $EX_1^2=\psi''(0+)<\i$, then $\lambda_\mu<\i$
if and only if $\int_0^\i x\mu(dx)<\i$ in which case $\int_{\R}|x|\mu(\td x)<\i$.
\end{prop}

The rest of this section is devoted to the proof of Proposition \ref{prop:min}.
The proof is based on two auxiliary results which are of independent interest:

\begin{lemma}\label{lem:Winf}
(i) If $X$ drifts to $+\i$ then $W$ is bounded by $1/\psi'(0+)<\i$.

(ii) If $X$ oscillates then $W(a)/a\to 2/\psi''(0+)$ as $a\to\i$.

(iii) Let $\L$ the L\'{e}vy measure of $X$ and
$C = 1 + \ovl W(a)\int_{-\i}^{-1}|x|\L(\td x)$. Then, if $X$ oscillates,
for any $a,x>0$,
\begin{equation}\label{eq:Winf}
a\wedge x-C\leq a\le(1 - \dfrac{W(a-x)}{W(a)}\ri)  \leq  x.
\end{equation}
\end{lemma}

\proof (i) If $X$ drifts to infinity, then $\psi'(0+)>0$ and
it follows by a Tauberian theorem (e.g.\ \cite[p.10]{bert})
and the definition of $W$ that $\lim_{x\to\i}W(x)=1/\psi'(0+)$.

(ii) Again appealing to a Tauberian theorem and noting that 
$\psi'(0+)=0$ in this case,
it follows that $W(x)/x\to 2/\psi''(0+)$ as $x\to \i$.

(iii) Since $E[X_1^+] \leq E[\sup_{t\leq 1} X_t] < \i$
(Bertoin \cite[Ch VII.1]{bert}) and $E[X_1^+] = E[X_1^-]$
(as $X$ oscillates), it follows that $E[|X_1|]$ is finite.

Write $\rho_{-a,0} = H_{\R\backslash[-a,0]}$ for the first exit time from $[-a,0]$.
We next show that the expectation of
the undershoot $o_a:=X_{\rho_{-a,0}} - X_{H_{-a,0}}$ is bounded
by $C$, which is finite since $E[|X_1|]<\i$. We have
\begin{align*}
E_{-x}\big[|o_a|\big] &\leq 1+E_{-x}\big[|o_a|
\mathbf{1}_{\{|o_a|\geq 1\}}\big]\leq 1+ E_{-x}\le[\sum_{s\ge 0}|\Delta X_s|
\mathbf 1_{\{\Delta X_s\leq -1, s=\rho_{-a,0}\}}\ri]\\
&\leq 1+ E_{-x}\le[\sum_{s\geq 0}|\Delta X_s|
\mathbf{1}_{\{\Delta X_s\leq -1\}}\mathbf{1}_{\{s\leq
\rho_{-a,0}\}}\ri]\\
&= 1 +  c \int_{-\i}^{-1}|z|\Lambda(\td z)\leq C,
\end{align*}
with $c = E_{-x}[\rho_{-a,0}]$, where we applied the compensation formula to the Poisson
point process $(\D X_s)$ and used that $c\leq \ovl W(a)$ 
(see \eqref{twosidedexp}). 

Using \eqref{twosided} yields that
\begin{equation}
E_{-x}[X_{\rho_{-a,0}}  - o_a]
= -a P_{-x}[\rho_{-a,0} < H_0]
= -a\le(1 - \frac{W(a-x)}{W(a)}\ri).
\label{eq:XHaa}
\end{equation}
Since $X$ oscillates and $E[|X_1|]<\i$, we note that $X$ is a martingale.
An application of Doob's optional
stopping theorem in conjunction with the dominated convergence
theorem implies that
\begin{equation*}
E_{-x}[X_{\rho_{-a,0} }] = \lim_{t\to\i}
E_{-x}[X_{t\wedge \rho_{-a,0} }] = - x.
\end{equation*}
The assertion in \eqref{eq:Winf} follows instantly from
\eqref{eq:XHaa}. Note in particular that for $x>a$ it trivializes.\exit\smallskip

\bigskip
\bigskip

We consider next the expected discounted local time up to $H_0$:
$$
v^q(x,y) = E_x\le[\I_0^{H_0}\te{-qt}\td L^y_t\ri],\quad q>0.
$$
Note that, when $q\downarrow 0$, $v^q(x,y)$ converges to $v(x,y)=E_x
[L^y_{H_0}]$, by monotone convergence. Recall
that the Laplace transform of the first hitting time
$H_0$ is given by (Bertoin \cite[Thm II.19]{bert})
\begin{equation}\label{eq:LTHE}
%E[\te{-qH_{\eta}}\mbf 1_{\{H_\eta < \i\}}] = 
E_{x}[\te{-qH_{0}}\mbf 1_{\{H_0< \i\}}] = 
\frac{u^q(-x)}{u^q(0)},\quad q>0,
\end{equation}
which is equal to $\te{\F(q)x}$ for $x<0$. 

\begin{lemma}\label{lem:on_blj}
Let $q>0$ and $y\in \R$. The following hold true:
\BI
\im[(i)] The process
$
v(X_t,y)+L^y_t-\frac{u(y)}{u(0)}L^0_t,
$
with $\frac{u(y)}{u(0)} := \lim_{q\downarrow0}\frac{u^q(y)}{u^q(0)}$,
is a martingale.
\im[(ii)] For $q>0$ it holds that
\begin{equation}\label{eq:uqvq}
v^q(x,y) = u^q(y-x) - \frac{u^q(y)u^q(-x)}{u^q(0)}.
\end{equation}
As a consequence, if $X$ does not drift to $-\i$,
\begin{equation}\label{eq:genv}
v(x,y) = W(x) + W(-y) - W(x-y) - W(x)W(-y)\psi'(0+).
\end{equation}
\im[(iii)] If $X$ drifts to $+\i$,
$E[\t_{L(\i)-}]= \psi''(0+)/[\psi'(0+)]^2$.
\EI
\end{lemma}
\proof
(i) In view of the strong Markov property, it follows that
\begin{eqnarray*}
E_x\le[\int_0^\i\te{-q t}\td L^y_t\ri] &=& v^q(x,y) +
E_x[\te{-q H_0}\mbf 1_{\{H_0 < \i\}}]E_0\le[\int_0^\i\te{-qt}\td L^y_t\ri].
\end{eqnarray*}
Taking note of (\ref{eq:LTHE}), \eqref{eq:uqy} and the spatial homogeneity
of a L\'{e}vy process, we see that \eqref{eq:uqvq} holds.
Applying again the Markov property it follows that, for $q>0$,
the process $M^q=\{M^q_u,u\ge 0\}$ is a UI martingale:
\begin{eqnarray*}
M^q_u &=& E\le[\int_0^\i\te{-q t}\td L^y_t - \int_0^\i\te{-q t}\td L^0_t
\frac{u^q(y)}{u^q(0)} \Big|\mc F_u\ri]\\
&=& \int_0^u\te{-q t}\td L^y_t - \int_0^u\te{-q t}\td L^0_t\frac{u^q(y)}{u^q(0)}\\
&\phantom{=}& + \te{-qu}\le[u^q(y-X_u) - u^q(-X_u)\frac{u^q(y)}{u^q(0)}\ri]\\
&=& \int_0^u\te{-q t}\td L^y_t - \frac{u^q(y)}{u^q(0)}
\int_0^u\te{-q t}\td L^0_t + \te{-qu}v^q(X_u,y).
\end{eqnarray*}
Observe that $M^q_u \to M^0_u$ pathwise a.s.\ as $q\downarrow0$. Furthermore, $v^q(x,y)\leq v^q(y,y)\leq v(y,y)$ and the expectation of the other two integrals is bounded
by $E_x[L^y_t + L^0_t]$ which is finite (cf.\ Section \ref{sec:levy}).
It now follows by the dominated convergence
theorem that $M^0_t=L^y_t -\frac{u(y)}{u(0)}L^0_t+v(X_t,y)$ is a martingale.
Notice that, if $X$ oscillates, $u(y)/u(0)=1$ and the resulting
martingale coincides with the one used by Bertoin and
Le Jan \cite{blj}.

(ii) In view of \eqref{eq:uq} and \eqref{eq:uqvq} it follows that
$v^q(x,y) = W^{(q)}(-y)\te{\F(q)x} + W^{(q)}(x)\te{-\F(q)y} - W^{(q)}(x-y)
- \Phi'(q)^{-1}W^{(q)}(-y)W^{(q)}(x)$. Letting then $q\downarrow 0$ equation
\eqref{eq:genv} follows by continuity of $W^{(\cdot)}$ and since
$\lim_{q\downarrow0} \F'(q) = \psi'(0+)^{-1}$ in the case when $X$ 
does not drift to $-\i$ (see
Section \ref{sec:levy}).

(iii) If $X$ drifts to $+\i$, it holds that $L(\i)$ follows an exponential
distribution and $(\t_\ell, \ell< L(\i))$ has the same
law as a subordinator $\WT\t$ killed at an independent exponential time,
in view of \cite[Thm. IV.8]{bert}.
Further, it follows by a change of variables, \eqref{eq:uqy} and
\eqref{eq:uq} that
\begin{equation}\label{eq:fpq}
\F'(q) = \int_0^\i E[\te{-q\t_\ell}\mbf 1_{\{\ell<L(\i)\}}]\td \ell.
\end{equation}
From \eqref{eq:fpq}  we deduce that $L(\i)\sim\exp(1/\F'(0))$ and
that Laplace exponent of $\WT\t$, $\k(q)= - \log E[\te{-q\WT\t_1}]$,
is equal to $\k(q)=\F'(q)^{-1} - \F'(0+)^{-1}$. In particular,
\begin{equation}\label{eq:exptt}
E[\WT\t_t] = t\k'(0+) = \le.\frac{\td}{\td q}\ri|_{q=0+} \frac{t}{\Phi'(q)} =
\frac{\psi''(0+)}{\psi'(0+)}t,
\end{equation}
using that $\psi'(\F(q)) = [\F'(q)]^{-1}$ and
$\F(0)=0$ if $X$ drifts to $+\i$.
Finally, if $\eta(a)$ denotes an independent exponential random variable
with parameter $a = 1/\F'(0)$, it follows from \eqref{eq:exptt}
and the relation between $\t$ and $\WT\t$ described above that
$E[\t_{L(\i)-}] = E[\WT\t_{\eta(a)}] = \int_0^\i a\te{-a\ell}E[\WT\t_\ell]\td \ell =
\psi''(0+)/\psi'(0+)^2$.
\exit\smallskip\\
\proof{\it of Proposition \ref{prop:min}}. 
Using \eqref{eq:genv} for $y>0$ we have that $V_\mu(y)=\int_0^\i (W(x)-W(x-y))\mu(\td x)$ which increases to $\lambda^+_\mu := \int_0^\i W(x)\mu(\td x)$ as $y\nearrow \i$. It follows that $\lambda_\mu\geq
\lambda_\mu^+$.

(ii) Suppose first that
$X$ drifts to infinity. Then, in view of Lemma \ref{lem:Winf},
(\ref{eq:genv}) and (\ref{eq:VBLJ}), we see that $V_\m$ is
bounded. In consequence $\lambda_\mu<\i$.

(iii) Part (i) which we prove below and the fact that 
$E[L_{\tilde{T}_{\f_\pm}}]=\int_0^{\i}
W(x)\mu(\td x)$, which we shall derive in \eqref{eq:est1} below,
immediately yield that $\lambda_\mu = \int_0^{\i}
W(x)\mu(\td x)$. Suppose that $X$ oscillates.
In view of the asymptotics in Lemma \ref{lem:Winf}(ii) 
and the fact that $W'(0)=2/\s^2<\i$, it then follows that the continuous function
$s\mapsto W(s)/s$ attains on $[0,\i]$ a finite positive maximum $c_+$ 
and finite non-negative minimum $c_-$ (which is positive if $\psi''(0+)<\i$). 
As a consequence,
\begin{eqnarray}\label{est1}
c_-\int_0^\i x\mu(\td x) &\leq& \int_0^\i W(x)\mu(\td x) \leq c_+\int_0^\i x\mu(\td x),\\
c_-\int_{-\i}^0 |x|\mu(\td x) &\leq& \int_{-\i}^0 W(-x)\mu(\td x) \leq c_+\int_{-\i}^0 |x|\mu(\td x).
\label{est2}
\end{eqnarray}
Next, noting that $W'(a)\leq W'(0)$, which follows by the fact (e.g.~\cite[eq. (13)]{mp-sem05} and 
\cite[Ch. VI]{bert}) that 
$$P_x\le[X_{H_{(-\i,0]}} = 0\ri] = \dfrac{W'(x)}{W'(0)}\q\text{for $x>0$},
$$
we deduce from \eqref{eq:assum} that 
\begin{equation}\label{est3}
\int_{-\i}^0 W(-x)\mu(\td x) \leq \int_0^\i W(x)\mu(\td x).
\end{equation}
The statements in (iii) follows by combining \eqref{est1}, \eqref{est2} and \eqref{est3}.

(i) Recall from Lemma \ref{lem:on_blj} that
$N^y_t=v(X_t,y)+L^y_t-\frac{u(y)}{u(0)}L^0_t$ is a martingale with
$N^y_0=0$. Localizing and using monotone and dominated convergence
theorems (note that $v(x,y)\leq v(y,y)$) we obtain from Doob's
optional sampling theorem that
$\frac{u(y)}{u(0)}E[L_S^0]=V_\mu(y)+E[L^y_S]\geq V_\mu(y)$. Taking
the supremum over $y$, and noting that $\frac{u(y)}{u(0)}\leq 1$, we deduce that $E[L^0_S]\geq \lambda_\mu$.
\exit\smallskip

Taking into account \eqref{eq:genv}, the solution proposed
by Bertoin and Le Jan \cite{blj}
is explicit. The only practical drawback is that one has to
observe simultaneously the local
times at all levels which might be hard to implement. In contrast,
solutions presented in
Theorems \ref{thm:main} and \ref{thm:main0} only depend on
the local time at zero,
the sign of the present excursion and the present position of the process.
We stress however
that, unlike the solution of Bertoin and Le Jan, these results
only apply to restricted sets of target measures.
We give now some examples.

\subsection{Examples}\label{sec:exam}
Consider $X_t  = B_t + t - J_t$ where $B$ is a standard Brownian motion
and $J$ is a compound Poisson process that jumps at rate 1 with
exponentially distributed jump sizes. Then $X$ is a martingale, it oscillates
and its Laplace exponent is given by
$\psi(\th) = \th^2/2 + \th - \th/(1+\th)
= \frac{\th^2 (3+\th)}{2(1+\th)}.$
The scale function of $X$ thus reads as
\begin{equation}\label{eq:ex2scale}
W(x) = \frac{2x}{3} + \frac{4}{9}(1-\te{-3x}),\quad x\ge 0.
\end{equation}

For $a<0<b$ let us determine which measures $\mu$ on $\{a,b\}$,
$\mu=(1-p)\delta_{a}+p\delta_{b}$ can be embedded in $X$.
Such measures are covered by Theorem \ref{thm:main} and
invoking assumption \eqref{eq:assum} we see that it must hold that
$$p W(b)=2(1-p)\frac{W(-a)}{W'(-a)}\quad
\textrm{thus}\quad p=\frac{2W(-a)}{W(b)W'(-a)+2W(-a)}.$$
In particular, for $a=-b$ we would have $p=3/(4+2\te{-3b})$.
The measure $\mu$ is embedded with the stopping time
$$\tld T_{a,b}=\inf\{t: X_t=b\textrm{ or }X_t=a\textrm{ and }
X_{|_{(\t_{L_t-},t)}}<0\}.$$
It is possible to work out explicitly the stopping time of Bertoin and Le Jan in
\eqref{eq:TBLJ} for this special case. Using $W'(-a)\leq 2$, we obtain $\lambda_\mu=pW(b)=V_\mu(b)>V_\mu(a)=(1-p)W(-a)$
and
$$T^\mu_{BLJ}=\inf\left\{t: X_t=b\textrm{ or }L^{a}_t>
\frac{pW(b)-(1-p)W(-a)}{p(1-p)W(b)}L^0_t\right\}.$$
The difference with $\tld T_{a,b}$ is in the mechanism that determines
whether to stop or not when visiting $a$.

In parallel, extending the approach of Theorem \ref{thm:main0}
to atomic measures,
we are led to consider $H_{a,b}=\inf\{t: X_t\in \{a,b\}\}$. Note that the
measure $\nu\sim X_{H_{a,b}}$ places less mass in $b$ than $\mu$.
One can verify independently, using \eqref{twosided}
and \eqref{eq:ex2scale}, that indeed
$$\nu(\{b\})=\frac{W(-a)}{W(b-a)}<\frac{2W(-a)}{W(b)W'(-a)+2W(-a)}
=\mu(\{b\}).$$
Finally, it is instructive to verify that the measure $\rho\sim X_{H_{\R\setminus[a,b]}}$
does not verify assumption \eqref{eq:assum} and can not be treated
in our setup.
Note that $\rho$ is given by
$$\rho(\td x)=p_+\delta_{b}(\td x)+p_-\delta_{a}(\td x)+(1-p_+-p_-)\te{x-a}
\mathbf{1}_{x<a}\td x,$$
where $p_+$ and $p_-$ are given by \eqref{twosided} and
\begin{equation}\label{twosc}
P_x[X_{H_{\R \setminus [a,b]}} = a] = \dfrac{W'(b-a)}{W'(0)}
\le[\dfrac{W'(x-a)}{W'(b-a)} -  \dfrac{W(x-a)}{W(b-a)}\ri],
\end{equation}

respectively (the latter identity was shown in \cite[Prop. 1]{mp-sem05}).
One can also compute explicitly $V_\rho$ to find that it is a constant
$V_\rho\equiv p_+W(b)$ and thus the stopping time introduced by
Bertoin and Le Jan
\cite{blj} coincides with $H_{\R\setminus [a,b]}$.

\section{Excursion theoretical calculations}\label{sec:exc_comp}

In this section we derive a number of excursion measure identities
which are essential to obtain the formulae presented in Sections
\ref{sec:sol1}-\ref{sec:sol3}. We first set the notation and
briefly recall the main concepts of the excursion theory of a
spectrally negative L\'{e}vy process away  from zero, following
Bertoin \cite[Ch. IV]{bert}.

The excursion process $e = \{e_\ell, \ell\ge 0\}$ of $X$ away from 0
is defined as
$$
e_\ell = (X_u: \t_{\ell^-} \leq u < \t_\ell)\q\text{if $\t_{\ell^-} < \t_\ell$}
$$
and $e_\ell = \partial$ (a graveyard state) else.
The excursion process takes values in the space
$\mathcal E \cup \mathcal E^{(\i)} \cup \{\partial\}$
where
$$
\mathcal E = \{\epsilon\in D[0,\i): \exists\zeta
= \zeta(\epsilon)<\i: \epsilon(s)\neq 0, s\in (0,\zeta)\}
$$
is the space of excursions with finite lifetime and
$$
\mathcal E^{(\i)} = \{\epsilon\in D[0,\i): \epsilon(s) > 0, s > 0\}.
$$
are the excursions with infinite lifetime. Write $\z=\z(\e)$
for the lifetime of an excursion $\e\in\mathcal E\cup \mathcal E^{(\i)}$ and define the sign of an excursion via
\begin{equation}\label{eq:sgn_def}
\text{sgn}(\epsilon)=\lim_{s\searrow 0} \frac{\epsilon(s)}{|\epsilon(s)|}.
\end{equation}
According to the fundamental result of It\^{o} \cite{ito}, if $0$ is recurrent,
$e$ is a Poisson point process
with characteristic measure denoted by $n$; if $0$ is transient,
$\{e_\ell, \ell\leq L(\i)\}$ is a Poisson point process stopped
at its first entrance into $\mathcal E^{(\i)}$.

The process $X$ (with $X_0=0$) is said to creep downwards
across $x<0$ if $X_{H_{(-\i,x)}} = x$. According to Millar \cite{millar73}
a spectrally negative L\'{e}vy process can creep
downwards {\em only} if its Gaussian coefficient $\s^2$ is not zero.
Therefore, if $\sigma = 0$, the L\'{e}vy process always enters $(-\i,0]$
by a jump and never enters $(-\i,0]$ by hitting $0$.
In this case an excursion away from 0 is thus {\it either} infinite and
all the time strictly positive (after the starting point) {\it or}
it is first positive and then jumps into $(-\i,0)$ and finally
returns to 0 (recalling that we assume that $X$ oscillates
or drifts to $+\i$). Note that according to our definition \eqref{eq:sgn_def}
the latter is also of positive sign. In the case that $\sigma>0$, there are two
additional forms of excursions, namely those that stay positive or negative
all the time and hit zero in finite time.

\subsection{Supremum and infimum}
In the literature there are several results on the law
of the supremum of excursions of spectrally negative L\'{e}vy processes:
Bertoin \cite{bertoin2XI} and Avram et al.\ \cite{AKP} calculated this law
for excursion of $X$ away from its infimum and away from its supremum respectively.
Lambert \cite{lambert} calculated the law, under the excursion measure,
of the supremum of the excursions away from a point of a spectrally
negative L\'{e}vy process confined to a finite interval.
In this section we extend these results by deriving
the joint law of the supremum $\maxe$ and the infimum $\mine$ of an excursion $\e$
away from zero,
$$
\maxe = \sup_{s\leq \zeta} \e_s\quad\quad
\mine = \inf_{s\leq \zeta} \e_s,
$$
under the excursion measure $n$.
\begin{lemma}
For $\delta,\eta>0$ it holds that
\begin{equation}\label{eq:gen_exc_law}
n(1-\mbf 1_{\{\maxe<\eta, \mine>-\d\}}\te{-q\zeta})
=\frac{W^{(q)}(\eta + \delta)}{W^{(q)}(\eta)W^{(q)}(\delta)}.
\end{equation}
In particular,
\begin{eqnarray}\label{eq:gen_exc_law_up}
n(\maxe\ge \eta) &=& \frac{1}{W(\eta)}, \\
n(\mine \leq -\delta) &=& \frac{1}{W(\delta)} - \psi'(0+)
\label{eq:gen_exc_law_ddown}\\
\label{eq:gen_exc_law_down} n(\maxe<\eta, \mine\leq -\d) &=&
\frac{1}{W(\eta)} \le[\frac{W(\eta+\d)}{W(\delta)} - 1\ri].
\end{eqnarray}
\end{lemma}
\proof Let us first derive the identities
\eqref{eq:gen_exc_law_up} --- \eqref{eq:gen_exc_law_down} as
consequences of \eqref{eq:gen_exc_law}.

Letting $q\downarrow 0$ in \eqref{eq:gen_exc_law} we see that
\begin{equation}\label{eq:gen_exc_int}
n(\maxe\ge\eta\ \text{or}\
\mine\leq-\delta \ \text{or}\ \zeta = \i) =
\frac{W(\eta+\d)}{W(\eta)W(\delta)}.
\end{equation}
As $X$ does not drift to $-\i$, either $X$ oscillates in which
case $\z < \i$ $n$-a.s., or $X$ drifts to infinity in which case
$\maxe\ge\eta$ - thus the previous display is also equal to
$n(\maxe\ge\eta\ \text{or}\ \mine\leq-\delta)$. The identity
\eqref{eq:gen_exc_law_up} then follows by letting $\d\to+\i$ in
\eqref{eq:gen_exc_int} and using that $W(\delta+\eta)/W(\delta)$
converges to 1 as $\delta\to\i$ if $X$ does not drift to $-\i$ (to
see the latter recall that $\lim_{x\to\i}W(x)=1/\psi'(0+)<\i$ if
$X$ drifts to infinity and refer to Lemma \ref{lem:Winf} if $X$
oscillates). To obtain \eqref{eq:gen_exc_law_down} we note that
$n(\maxe<\eta, \mine\leq -\d) = n([\maxe<\eta, \mine>-\d]^c) -
n(\maxe\geq \eta)$ and employ \eqref{eq:gen_exc_int} and
\eqref{eq:gen_exc_law_up}. Note that \eqref{eq:gen_exc_law_ddown}
follows by letting $\eta\to\i$ in \eqref{eq:gen_exc_law_down}.

To show \eqref{eq:gen_exc_law} the first step is to establish
a link between the excursion measure and the expected discounted local
time. We show that for $\eta,\delta>0$
\begin{equation}\label{eq:exc_H_id}
n(1-\mbf 1_{\{\maxe<\eta,\mine>-\delta,\z<\i\}}\te{-q\zeta}) =
\le[E\le[\I_0^{H_{\eta,-\delta}}\te{-q t} \td L_t\ri]\ri]^{-1}.
\end{equation}
Indeed, changing from time scale  it follows that
\begin{multline*}
\lefteqn{E\le[\I_0^{H_{\eta,-\delta}}\te{-q t} \td L_t\ri] =
E\le[\int_0^\i \te{-q\t_\ell}\mbf 1_{\{\t_\ell < H_{\eta,-\delta}\}}
\td \ell\ri]}\\
= \I_0^\i \td \ell E\le[\exp\le(-q\sum_{0\leq u\leq \ell}(\t_u -
\t_{u^-}) \chi_{\{\ovl{e_u}<\eta,\unl{e_u}>-\delta\}}\ri)\ri]
\end{multline*}
where $\chi_A(\w)$ is the indicator of the set $A$ that is one if
$\w\in A$ and $+\i$ else. By the exponential formula for Poisson
point processes it then follows that the last line of the previous
display is equal to
\begin{multline*}
\I_0^\i\td \ell \exp \le(-\ell\,
n(1-\te{-q\z\chi_{\{\maxe<\eta,\mine>-\delta,\z<\i\}}})\ri)\\
= \le[n(1-\te{-q\z}\mbf
1_{\{\maxe<\eta,\mine>-\delta,\z<\i\}})\ri]^{-1}.
\end{multline*}
where we used that $\te{-\i}=0$. By the Markov property
$$
E\le[\I_0^\i \te{-q u}\td L_u\Big|\mc F_t\ri] = \I_0^t\te{-qu}\td L_u
+ \te{-qt} u^q(-X_t).
$$
Applying the optional stopping theorem at $H_{\eta,\delta}$ shows
that
\begin{eqnarray*}
E\le[\I_0^{H_{\eta,-\delta}}\te{-qt}\td L_t\ri] &=& u^q(0) -
u^q(-\eta) E[\te{-qH_{\eta,-\delta}}
\mbf 1_{\{H_\eta < H_{-\delta}, H_{\eta,-\delta}<\i\}}]\\
&&+\; u^q(+\delta)E[\te{-qH_{\eta,-\delta}}\mbf 1_{\{H_\eta >
H_{-\delta}, H_{\eta,-\delta}<\i \}}].
\end{eqnarray*}
Inserting now the expression of the Laplace transform of the
first hitting time $H_{\eta,-\delta}$ in terms
of the potential density $u^q$ from Pistorius \cite[Cor. 3]{mp-sem05},
it follows that
\begin{eqnarray*}
\lefteqn{E\le[\I_0^{H_{\eta,-\delta}}\te{-qt}\td L_t\ri]}\phantom{isssssss}\qquad\\
&=& u^q(0) -
\frac{u^q(-\eta)[u^q(\eta)u^q(0) - u^q(-\d)u^q(\d+\eta)]}%
{u^q(0)^2 - u^q(\d+\eta)u^q(-\delta-\eta)}\\
&&-\; \frac{u^q(\d)[u^q(-\d)u^q(0) - u^q(\eta)u^q(-\d-\eta)]}%
{u^q(0)^2 - u^q(\d+\eta)u^q(-\delta-\eta)}
\end{eqnarray*}
In view of \eqref{eq:uq} the identity \eqref{eq:gen_exc_law}
follows after some algebra.\exit

Next we turn to the calculation of the laws of $H_\eta$ and
$H_{-\delta}$ under the excursion measure $n$. When excursion's time scale is
considered $H_\eta$ refers to $H_\eta(\e)=\inf\{s: \e(s)=\eta\}$.

\begin{lemma}\label{lem:n_mm}
For $\eta,\delta>0$ it holds that
\begin{eqnarray}
\label{eq:t-un}
n(\te{-q H_{\eta}} \mbf 1_{\{\ovl\e \ge \eta\}}) &=& 1/W^{(q)}(\eta)\\
\label{eq:t-dn} n(\te{-q H_{-\d}} \mbf 1_{\{\maxe < \eta, \mine \leq
-\delta\}})
&=& \frac{\te{\F(q)\d}}{W^{(q)}(\eta)}\le[\frac{W^{(q)}(\d + \eta)}%
{W^{(q)}(\delta)} - \te{\F(q)\eta}\ri]
\end{eqnarray}
In particular, writing $W\star W$ for the convolution,
\begin{eqnarray}
\label{eq:exp_H+} n(H_{\eta}\mbf 1_{\{\maxe \ge \eta\}})
&=& \frac{W\star W(\eta)}{W(\eta)^2}\\
\nn n(H_{-\d} \mbf 1_{\{\maxe < \eta, \mine \leq -\delta\}}) &=&
\le[\frac{W\star W(\eta)}{W(\eta)^2} - \frac{\d\F'(0)}{W(\eta)}\ri]
\le[\frac{W(\eta+\delta)}{W(\d)} - 1\ri]\\
- \frac{1}{W(\eta)}&&\hspace*{-1.0cm}\le[\frac{W\star W(\eta+\d)}{W(\d)} -
\frac{W(\eta+\delta)W\star W(\d)}{W(\d)^2} - \F'(0)\eta\ri] \label{eq:exp_H-}\quad
\end{eqnarray}
\end{lemma}
\proof The expressions \eqref{eq:exp_H+} and \eqref{eq:exp_H-}
follow by evaluating the derivative with respect to $q$ at $q=0$
of the expressions \eqref{eq:t-un} and \eqref{eq:t-dn} respectively,
using the series representation (\ref{eq:expW}) of $W^{(q)}$.

Appealing to the compensation formula it follows that for $\eta>0$
$$
E[\te{-q H_\eta}] = E\le[\I_0^{H_\eta}\te{-q t}\td L_t\ri]
n(\te{-q H_{\eta}(\e)}\mbf 1_{\{\ovl\e \ge \eta\}}).
$$
In view of \eqref{eq:LTHE}, \eqref{eq:gen_exc_law} and \eqref{eq:exc_H_id}
it follows that the left-hand side is equal to $\te{-\F(q)\eta}$
and the first factor on the right-hand side is equal to
$\te{-\F(q)\eta}W^{(q)}(\eta)$. As a consequence, we see that
the expression \eqref{eq:t-un} holds true.

Similarly, an application of the compensation formula shows that
\begin{multline*}
E[\te{-q H_{\eta, -\delta}}] = E\le[\I_0^{H_{\eta,-\delta}}\te{-q
t}\td L_t\ri]\\ \times [n(\te{-q H_{\eta}(\e)} \mbf 1_{\{\ovl\e
\ge \eta\}}) + n(\te{-q H_{-\d}(\e)} \mbf 1_{\{\ovl\e < \eta,
\unl\e \leq -\delta\}})].
\end{multline*}
The first factor on the right-hand side is equal to
the reciprocal of \eqref{eq:gen_exc_law}, while
a short calculation employing \cite[Cor. 3]{mp-sem05}
and \eqref{eq:uq} shows that the left-hand side is equal to
$$
\frac{W^{(q)}(\eta+\delta)\te{\F(q)\d} +
W^{(q)}(\d)(1-\te{\F(q)(\d+\eta)})}{W^{(q)}(\d+\eta)}
$$
and the proof of \eqref{eq:t-dn} is complete.\exit

\subsection{Further computations in the presence of a Gaussian component}
Assume now specifically that $\s^2 > 0$: a Gaussian component is present.
In this case the process $X$ can creep both to positive
and negative levels and we can split excursions into negative
and positive ones. Recall that we defined the sign of an excursion $\e$ as
its sign at $t=0+$, i.e. $\text{sgn}(\e) =
\lim_{s\downarrow 0}\e(s)/|\e(s)|$. The signed maximum
functional then reads as
\begin{equation}\label{eq:def_sgnmax}
M(\e) = \frac{1}{2}(\text{sgn}(\e) + 1)\sup_{s\leq\zeta(\e)}\e(s)
+ \frac{1}{2}(\text{sgn}(\e) - 1)\inf_{s\leq\zeta(\e)}\e(s).
\end{equation}
In a positive excursion we thus only look at the maximum (and ignore
the infimum that may be attained in the excursion) and find
the following results for the law of $M$ under $n$:
\begin{lemma} Assume $\sigma>0$. For $a>0$, it holds that
\begin{equation}\label{eq:excmres}
n( M > a) = \frac{1}{W(a)}\q n(M < -a) =
\frac{W'(a)}{W(a)} \cdot \frac{1}{W'(0)}.
\end{equation}
\end{lemma}
\proof The first identity in \eqref{eq:excmres} follows from \eqref{eq:t-un}.
Since in a {\it negative} excursion $\maxe=0$, the second identity
in \eqref{eq:excmres} follows  by taking the limit $\eta\downarrow 0$
in \eqref{eq:gen_exc_law_down}
\begin{eqnarray}
n(\maxe=0, \mine\leq -\d) = \lim_{\eta\downarrow 0} n(\maxe<\eta,
\mine\leq -\d) = \frac{1}{W'_{}(0)}
\frac{W'_{}(\delta+)}{W_{}(\delta)} \label{eq:exc_e0_lawdown}
\end{eqnarray}
\exit
\\Finally, we record for later use  the form of
the first moment of $H_{-\d}$ under $n$:
\begin{lemma}\label{lem:nm_mm}
Assume $\sigma>0$ and $\delta,\eta>0$. It holds that
\begin{multline*}
n(H_{-\delta}\mbf 1_{\{\ovl\e = 0, \unl\e \leq -\delta\}})\\
= \frac{1}{W'(0)}\le[\F'(0)\le(1 - \d\frac{W'(\d)}{W(\d)}\ri)
+ \frac{W'(\d)(W\star W)(\d)}{W(\d)^2} - \frac{(W\star W')(\d)}{W(\d)}\ri].
\end{multline*}
\end{lemma}
\proof Rewriting (\ref{eq:t-dn}) as 
$$
\frac{\eta\te{\Phi(q)\delta}}{W^{(q)}(\eta)} 
\le[\frac{W^{(q)}(\d + \eta) - W^{(q)}(\d)}{\eta W^{(q)}(\d)} + \frac{1 - \te{\Phi(q)\eta}}{\eta}\ri]
$$
and then taking the limit $\eta\downarrow 0$ shows that
$$
n(\te{-q H_{-\d}} 1_{\{\ovl\e = 0, \unl\e \leq -\delta\}})
= \frac{\te{\F(q)\d}}{W^{(q)\prime}(0)}\le[\frac{W^{(q)\prime}(\d)}%
{W^{(q)}(\delta)} - \F(q)\ri].
$$
The identity follows by subsequently calculating the right-derivative
with respect to $q$ in $q=0$, using the series representation \eqref{eq:expW}
of $W^{(q)}$ and the facts that $\F(0)=0$ if $X$ does not drift to $-\i$ 
and that $W^{(q)\prime}(0)=W'(0)$, again in view of \eqref{eq:expW}.
\exit

\section{Proofs of Theorems 1, 2 and 3}\label{sec:proofs}

\proof \emph{of Theorem \ref{thm:main}}

As Gaussian component is present ($\sigma>0$) the process $X$ can creep both to positive
and negative levels and excursions are either positive
and hit zero, negative and hit zero, positive and jump below
zero and then hit zero. Recall our definition of the sign of an excursion given in \eqref{eq:sgn_def} and the signed maximum
functional in \eqref{eq:def_sgnmax} and note that we only look at the infimum along negative excursions. The process $M(e_\ell)$ is a Poisson point process with characteristic measure $n(M\in
\td a )$.
The embedding part of the theorem can be proved similarly as Theorem 1 in Ob\l\'oj \cite{ob_gen} (accounting for atoms as in Cox, Hobson and Ob\l\'oj \cite{cho}) using the excursion measure calculations in \eqref{eq:excmres}.
Before carrying out the proof let us specialize briefly to the case of non-atomic measures. Define for $x<0<y$
$$
D_\m(y)=\I_{[0,y]}W(s)\m(\td s) \q\text{and}\q
G_\m(x)=W'(0)\I_{[x,0]}\frac{W(-s)\m(\td s)}{W'(-s)}\, .$$ Then $\tld\psi(x)$ is given by $\psi_+(x)$ for $x\ge 0$ and by $-\psi_-(x)$ for $x<0$ with
\begin{equation}\label{eq:fdef0atoms}
\begin{split}
\psi_+(y) =& \I_0^y \frac{W(s)\mu(\td s)}{(1+\ovl\mu(s) - \ovl\mu(g(s)))}\\
\psi_-(z) =& \I_{-z}^0 \frac{W'(0)W(-s)\mu(\td s)}{W'(-s) (1+\ovl\mu(f(s)) - \ovl\mu(s))}
\end{split} \end{equation}
where $g(s) = G^{-1}_\m(D_\m(s))$ and $f(s) = D^{-1}_\m(G_\m(s))$. The assumption \eqref{eq:assum} simplifies to $D_\m(\i) = G_\m(-\i)$ which is the admissibility criterion of Ob\l\'oj \cite{ob_gen}.\\
For the remainder of the proof, we write $T=\tld T_{\f_\pm}$. Using Poisson Point Process properties of the excursion process, we have that
\begin{equation}
\label{eq:lawLT}
P(L_T > k) = \exp\le(-\I_0^kn(M>\f_+(\ell))+n(M<-\f_-(\ell))\td \ell\ri).
\end{equation}
Recall the definitions of $\a, \x,\tld \psi$ and the excursion measure calculation given in \eqref{eq:excmres}.
It follows that, for $a_*<a<1$, $\f_+(\x(a))=F^{-1}_\mu(a)$
and $\f_-(\x(\a(a)))=-F^{-1}_\mu(\a(a))$. Finally note that
$$\a'(a)=-\frac{W'(-F^{-1}_\mu(\a(a)))W(F^{-1}_\mu(a))}{W'(0)W(-F^{-1}_\mu(\a(a)))}
=-\frac{n(M<F^{-1}_\mu(\a(a)))}{n(M>F^{-1}_\mu(a))}.$$
In consequence we get
\begin{eqnarray}
P(L_T>\x(a))&=&\exp\left(-\int_{a_*}^a\left(\frac{n(M<F^{-1}_\mu(\a(u)))}{n(M>F^{-1}_\mu(u))}+1\right)\frac{\td
u}{1-u+\a(u)}\right)\nonumber\\
&=&\exp\left( \ln[1-u+\a(u)]\Big|_{a_*}^a\right)=1-a+\a(a).\label{eq:lawTinH}
\end{eqnarray}
Let $a_*<c_0<1$ be such that $n(M<F^{-1}_\mu(\a(c_0)))+n(M>F^{-1}_\mu(c_0))\leq 1$. Then we can write
\begin{eqnarray*}
\x(1)&=&2\int_{a_*}^1\frac{du}{n(M>F^{-1}_\mu(u))(1-u+\a(u))}\\
&\geq& \int_{c_0}^1\frac{n(M<F^{-1}_\mu(\a(u)))+n(M>F^{-1}_\mu(u))}{n(M>F^{-1}_\mu(u))(1-u+\a(u))}\td
u\\
&=&-\ln(1-u+\a(u))|_{c_0}^1=\i.
\end{eqnarray*}
Likewise we show that $\x(0)=-\i$. We conclude via \eqref{eq:lawTinH} that $L_T<\i$ a.s. This readily implies that $T<\i$ a.s.\ when $X$ oscillates. If $X$ drifts to $+\i$ then it a.s.\ hits the level $\f_+(L_\i)<b_\mu$ and thus $T<\i$ a.s.\\
Let $x>0$. Then we can write
\begin{eqnarray*}
P(X_T>x)&=&\int_0^\i \P(L_T>l)n(M>\f_+(l))\mathbf{1}_{\f_+(l)>x}\td l\\
&=&\int_{a_*}^1 P(L_T>\x(a))n(M>F_\mu^{-1}(a))\mathbf{1}_{F_\mu^{-1}(a)>x}\td \x(a)=\mu((x,\i))
\end{eqnarray*}
and likewise for $x<0$
\begin{eqnarray*}
P(X_T<x)&=&\int_0^\i \P(L_T>l)n(M<-\f_-(l))\mathbf{1}_{-\f_-(l)<x}\td l\\
&=&\int_{a_*}^1 P(L_T>\x(a))n(M<F_\mu^{-1}(\a(a)))\mathbf{1}_{F_\mu^{-1}(\a(a))<x}\td \x(a)\\
&=&\int_{a_*}^1
\frac{n(M<F_\mu^{-1}(\a(a)))}{n(M>F^{-1}_\mu(a))}\mathbf{1}_{F_\mu^{-1}(\a(a))<x}\td
a=\mu(-\i,x),
\end{eqnarray*}
which proves that $X_T\sim\mu$.

We turn to the proof of the minimality of the stopping time $T$. Let $S\leq
T$ be a stopping time with $X_S\sim X_T$. We will show that $S=T$ a.s. We start by computing
$E[L_T]$:
\begin{eqnarray}
\nn
E[L_T] &=& \I_0^\i P(L_T\ge k)\td k = \int_{a_*}^1 P(L_T>\xi(u))d\xi(u)\\
%\I_{\f_+(0)}^{\f_+(\i)}P(L_T\ge \psi_+(k))\td\psi_+(k)\\
&=& \int_{a_*}^1\frac{du}{n(M>F^{-1}_\mu(u))}=\I_0^\i W(y)\m(\td y).
\label{eq:est1}
\end{eqnarray}
From Proposition \ref{prop:min} we deduce that $E[L_S]\geq E[L_T]$ and since
$0\leq L_S\leq L_T$ we conclude that $L_S=L_T$ a.s.\ that is $S$ and $T$
happen in the same excursion away from zero. From the definition of $T$
and since $S\leq T$ we see that $\text{sgn}(X_S)=\text{sgn}(X_T)$.
Absence of positive
jumps implies $X_S\mathbf{1}_{X_S\geq 0}\leq X_{T}\mathbf{1}_{X_{T}\geq 0}$ a.s.\
and in consequence $S=T$ on the set $\{X_S\geq 0\}=\{X_{T}\geq 0\}$. For the
negative values we have to deal with the undershoot.
Let $\varrho=\inf\{t:X_t\leq \varphi_-(L_t)\textrm{ and  sgn}(e_{L_t})=-1\}$.
Note that on $\{X_T<0\}=\{\t_{T-}<\varrho\leq T\}$ we have $X_u<X_T$ for $u\in (\varrho,T)$.
If $P(S\leq \varrho\leq T)=0$ then $\varrho<S\leq T$ on $\{X_T<0\}$. Thus
$X_S\leq X_T$ on $\{X_T<0\}=\{X_S< 0\}$ and
since $X_S\sim X_T$ we deduce that $S=T$.
Suppose next that $P(S\leq \varrho\leq T)=\e>0$. Then,
working conditionally on $\{S\leq \varrho\leq T\}$, we apply the Markov property
at $S$ to see that, starting from $X_S$, there is a positive probability of
hitting zero before hitting $[-\f_-(L_S),-\i)$ which in turn means
that $P(L_T>L_S)>0$ which gives the contradiction. We conclude that
$S=T$ a.s. and therefore $T$ is minimal.

Next we show the finiteness of $E[T]$. Recall that we assume that $X$ drifts
to $+\i$. Conditioning on $L_T$ we can write $E[T]$ as
\begin{eqnarray*}
E[T] &=& E[E[T|L_T]]\\
&=& E[\t_{L_T -}] + E\Big[E[(T-\t_{L_T-})\mbf 1_{\{X_T=\f_+(L_T)\}}|L_T]\\
&&+\; E[(T-\t_{L_T -})\mbf 1_{\{X_T=-\f_-(L_T)\}}|L_T]\Big].
\end{eqnarray*}
Properties of Poisson point processes imply that
$$p_+(k):=P(X_T=\f_+(k)|L_T=k)=\frac{n(\ovl \e >\f_+(k))}{n(\ovl \e >\f_+(k))+n(\unl\e \leq -\f_-(k))},$$
with an analogous expression for $p_-(k)=1-p_+(k)$.
Since the law of the first (positive) excursion away from zero
with supremum larger than $\f_+(k)$ is given by
$n(\cdot\ |\ovl\e > \f_+(k))$ and, conditional on $L_T$,
our functional of the first excursion is independent of $L_T$
(and similarly for the first negative excursion with infimum
smaller than $-\f_-(k)$) it holds that
\begin{eqnarray}
\nn E[T] &=& E[\t_{L_T -}] + \I_0^\i P(L_T\in\td k) \Big[n(H_{\f_+(k)}(\e)|\ovl\e \ge \f_+(k))
p_+(k) \\
\nn &&+\; n(H_{-\f_-(k)}(\e)|\ovl\e=0, \unl\e \leq - \f_-(k)) p_-(k)\Big]\\
\nn &=& E[\t_{L_T - }] + \I_0^\i P(L_T\ge k)\Big[
n(H_{\f_+(k)}(\e)\mbf 1_{\{\ovl\e \ge \f_+(k)\}})\\
&&+\; n(H_{-\f_-(k)}(\e)\mbf 1_{\{\ovl\e=0, \unl\e \leq - \f_-(k)\}})\Big]\td k
\label{eq:est0}
\end{eqnarray}
where in the last line used the form of $P(L_T\in\td k)$
that was displayed in \eqref{eq:lawLT}.

To show that $E[T]$ is finite we continue now by
estimating the three terms in the above display.
For the first term note that $E[\t_{L_{T-}}] \leq E[\t_{L(\i)-}]$
which is finite if $E[X_1^2]<\i$, in view of
Lemma \ref{lem:on_blj}(iii).
Changing variables in a similar way as in \eqref{eq:est1} yields that
\begin{equation}\label{eq:est2}
\I_0^\i P(L_T\ge k) \f_+(k)\td k = \I_0^\i yW(y)\m(\td y)
\end{equation}
and
\begin{equation}\label{eq:est3}
\I_0^\i P(L_T\ge k) \f_-(k)\td k = W'(0)\I_{-\i}^0 y\frac{W(-y)}{W'(-y)}\m(\td y).
\end{equation}
Further, noting that $W\star W (x) \leq x W(x)^2$ and
$(W\star W')(x)\leq W(x)^2$ (using that $W$ is increasing)
the statement of the Theorem now follows given the explicit forms of
$n(H_{\f_+(k)}(\e)\mbf 1_{\{\ovl\e \ge \f_+(k)\}})$
and $n(H_{-\f_-(k)}(\e)\mbf 1_{\{\ovl\e=0, \unl\e \leq - \f_-(k)\}})]$
derived in Lemmas \ref{lem:n_mm} and \ref{lem:nm_mm}.

\exit\medskip\\
\proof \emph{of Theorem \ref{thm:main0}}

There are several ways in which we can stop. Firstly, an
excursion can start positive and have a maximum larger
than $\f_+(L_T)$. Secondly, an excursion can start positive,
have a maximum smaller than $\f_+(L_T)$ then jump negative and
have an infimum smaller than $-\f_-(L_T)$. Finally, an excursion
can start negative to achieve an infimum smaller than $-\f_-(L_T)$.
The last scenario is possible iff $\s>0$.

From standard considerations we obtain the law of $L_T$:
$$
P(L_T > k) = \exp\le(-\I_0^kn(\maxe \ge \f_+(s)\ \text{or}\
\mine\leq-\f_-(s))\td s\ri).
$$
Given $L_T=k$ either $X_T = \f_+(k)$ or $X_T=-\f_-(k)$. By the
property of Poisson point processes it thus follows that
\begin{eqnarray*}
p_+(k) &=& \frac{n(\maxe\ge \f_+(k))}{n(\maxe \ge \f_+(k)\ \text{or}\
\mine\leq-\f_-(k))}
\end{eqnarray*}
and $p_-(k) = P(X_T = -\f_-(L_T)|L_T=k) = 1 - p_+(k)$. Also, for $h:\R \to \R_+$,
$$
E[h(X_T)] = \I_0^\i P(L_T\in\td k)\left[ h(-\f_-(k))p_-(k) + h(\f_+(k))p_+(k)\right]\, .
$$
By choosing $h(z)=\mbf 1_{z\ge y}$ for $y>0$ and writing
$\psi_\pm$ for the inverses of $\f_\pm$ we get
$$
\ovl\mu(y) = \I_{\psi_+(y)}^\i \td k\ n(\maxe\ge \f_+(k))P(L_T\ge k)
$$
and (with $h(z)=\mbf 1_{z\ge x}$ ($x<0$))
$$
\ovl\mu(x) = P(L_T\leq \psi_-(-x)) + \I_{\psi_-(-x)}^\i \td k\
n(\maxe\ge \f_+(k))P(L_T\ge k).
$$
Reasoning as in the proof of Thm.~1 in Ob{\l}\'{o}j \cite{ob_gen}
we find that
$$
\td \psi_+(y) = \frac{-\td\ovl\mu(y)}{n(\maxe\ge y)(1 + \ovl\mu(y) - \ovl\mu(g(y)))},
$$
where $g(y) = -\f_-(\psi_+(y))$ and
$$
\td \psi_+(y) = \frac{\td\ovl\mu(g(y)) - \td\ovl\mu(y)}%
{n(\maxe\ge y\ \text{or}\ \mine\leq g(y))(1 + \ovl\mu(y) - \ovl\mu(g(y)))}.
$$
Comparing these two expressions shows that
$$
\frac{\td\ovl\mu(y)}{n(\maxe\ge y)} =
- \frac{\td\ovl\mu(g(y))}{n(\maxe<y, \mine\leq g(y))}.
$$
This leads to the following equation for $g$ that must be satisfied:
$$
\ovl\mu(g(x))-\ovl\mu(0) =
\I_0^x \frac{n(\maxe<y, \mine\leq g(y))}{n(\maxe\ge y)}\mu(\td y).
$$

Assuming that $\m$ is absolutely continuous w.r.t.
the Lebesgue measure and writing $f_\mu(x)$ for its density at $x$ it follows from \eqref{eq:gen_exc_law_up} and \eqref{eq:gen_exc_law_down} that $g:\R_+\to\R_-$ must satisfy \eqref{eq:g_diff}
%\begin{equation}\label{eq:g_diff_copy}
%\frac{\td g}{\td y}(y) = - \frac{W(y - g(y)) - W(-g(y))}{W(-g(y))}
%\frac{f_\mu(y)}{f_\mu(g(y))}\end{equation}
with $g(0) = - \f_-(\psi_+(0)) = 0$. We see that if such $g$ exists then it is plainly a decreasing
function, as required. Furthermore, for $g$ and its inverse to be well defined we have to have
$g(x)\to a_\mu$ as $x\to b_\mu$, which is the analogue of criterion \eqref{eq:assum} in Theorem \ref{thm:main}.
The formulae in Theorem \ref{thm:main0} then follow.\\
It remains to see that $T$ is minimal. Let $S\leq T$ with $X_S\sim X_T$. Reasoning presented in the proof of Theorem \ref{thm:main} applies if we can show that $L_S=L_T$ and $sgn(X_S)=sgn(X_T)$. To this end, note that $P (L_T>\psi_+(y))=1+\ovl \mu(y)-\ovl\mu(g(y))$, so that
$$E[L_T]=\int_0^\i P(L_T>k)=\int_0^\i P(L_T>\psi_+(k))d\psi_+(k)=\int_0^\i W(y)\mu(\td y)$$
which by Proposition \ref{prop:min} is the lower bound on $E L_S$. We thus have $0\leq L_S\leq L_T$ with $E[L_S]=E[L_T]$ and thus $L_S=L_T$ a.s. From the definition of $T$ we see promptly that $\{X_S\geq 0\}=\{X_T\geq 0\}$,
 and $X_S\sim X_T$ implies $sgn(X_S)=sgn(X_T)$ a.s.
\exit\smallskip\\
\proof \emph{of Theorem \ref{thm:pos}}

As the running supremum of an excursion is continuous, the embedding
part of the theorem follows directly
from Ob\l\'oj \cite{ob_gen}. So does the minimality of $T_{\f_\mu}$ and
the statement $\sup_{t\leq T_{\mu}} X_t=X_{T_\mu}$ is immediate. \exit
\subsection*{Acknowledgement} 
We thank an anonymous referee for his careful reading of the paper and
helpful remarks.

\end{document}